\documentclass[11pt]{article}
\usepackage{mathpple}
\usepackage{amsfonts}
\usepackage{amsmath}
\usepackage{amstext}
\usepackage{amssymb}
\usepackage{amsthm}
\usepackage{amscd}
\usepackage{xypic}

\usepackage{vmargin}
\def\sp{{\mathop{\rm Sp}}}
\def\gsp{{\mathop{\rm GSp}}}
\def\orth{{\mathop{\rm Orth}}}
\def\O{{\mathop{\rm O}}}
\def\u{{\mathop{\bf U}}}

\def\qed{$\square$}

\global\let\bar\undefined
\newcommand{\bar}[1]{{\overline{#1}}}

\global\let\rank\undefined
\DeclareMathOperator{\rank}{rank}

\def\geom{{\rm geom}}
\def\cl{{\rm cl}}
\def\cgn{{_N\calc_g}}
\def\mgn{{_N\calm_g}}

\def\defeq{:=}

\def\ff{{\mathord{\mathbb F}}}

\renewcommand{\thesubsection}{\ifnum\value{subsection}=0{\arabic{section}}\fi\ifnum\value{subsection}>0{\arabic{section}.\arabic{subsection}}\fi}

\numberwithin{equation}{subsection}
\numberwithin{table}{section}

\newtheorem{theorem}[equation]{Theorem}
\newtheorem{lemma}[equation]{Lemma}
\newtheorem{prop}[equation]{Proposition}
\newtheorem{cor}[equation]{Corollary}

\theoremstyle{remark}
\newtheorem{situation}[equation]{Situation}

\def\mypara{\paragraph}

\def\ang{\ip}
\def\calc{{\cal C}}
\def\cald{{\cal D}}
\def\cale{{\cal E}}
\def\calf{{\cal F}}
\def\calh{{\cal H}}
\def\calm{{\cal M}}
\def\caln{{\cal N}}
\def\calo{{\cal O}}
\def\calp{{\cal P}}
\def\calt{{\cal T}}
\def\calu{{\cal U}}
\def\cross{\times}
\DeclareMathOperator{\aut}{Aut}
\DeclareMathOperator{\End}{End}
\DeclareMathOperator{\frob}{Fr}
\def\ff{{\mathord{\Bbb F}}}
\DeclareMathOperator{\gl}{GL}
\def\gth{\frak}
\DeclareMathOperator{\id}{id}
\def\inject{\hookrightarrow}
\def\integ{{\mathord{\Bbb Z}}}
\def\inv{^{-1}}
\def\iso{\cong}
\DeclareMathOperator{\jac}{Jac}
\DeclareMathOperator{\mat}{Mat}
\DeclareMathOperator{\mult}{mult}
\def\nat{{\mathord{\Bbb N}}}
\def\pf{\mypara{Proof}}
\def\ra{\rightarrow}
\def\rat{{\mathord{\Bbb Q}}}
\DeclareMathOperator{\spec}{Spec}
\def\units{^\cross}
\newcommand{\abs}[1]{{\left|#1\right|}}
\newcommand{\floor}[1]{{\lfloor #1 \rfloor}}
\newcommand{\half}[1]{\frac{#1}{2}}
\newcommand{\ip}[1]{{{\langle #1 \rangle}}}
\newcommand{\lie}[1]{{\gth #1}}
\newcommand{\oneover}[1]{\frac{1}{#1}}
\newcommand{\rest}[1]{|_{#1}}
\newcommand{\st}[1]{\{#1\}}
\newcommand{\til}[1]{{\widetilde{#1}}}

\title{The distribution of class groups of function fields}
\author{Jeffrey D. Achter\\
{\tt j.achter@colostate.edu}
}
\date{\empty}
\begin{document}
\maketitle

\begin{abstract}
For any sufficiently general family of curves over a finite field
$\ff_q$ and any elementary abelian $\ell$-group $H$ with $\ell$
relatively prime to $q$, we give an
explicit formula for the proportion of curves $C$ for which
$\jac(C)[\ell](\ff_q)\iso H$.  In doing so, we prove a
conjecture of Friedman and Washington.
\end{abstract}

%\tableofcontents

In 1983, Cohen and Lenstra introduced heuristics \cite{cl} to explain
statistical observations about class groups of imaginary quadratic
fields.  Their principle, although still unproven, remains an
important source of guidance in number theory.  A concrete application
of their heuristics predicts that an abelian group occurs as a class
group of an imaginary quadratic field with frequency inversely
proportional to the size of its automorphism group.

Six years later Friedman and Washington \cite{fw} addressed the
function field case.  Fix a finite field $\ff_q$ and an abelian
$\ell$-group $H$, where $\ell$ is an odd prime relatively prime to
$q$.  Friedman and Washington conjecture that $H$ occurs as the
$\ell$-Sylow part of the divisor class group of function fields over
$\ff_q$ with frequency inversely proportional to $\abs{\aut(H)}$.  As
evidence for this, they prove that the uniform distribution of
Frobenius automorphisms of curves of genus $g$ in $\gl_{2g}(\integ_\ell)$ would imply
their conjecture.  (Of course, autoduality of the Jacobian means that
these Frobenius elements are actually in $\gsp_{2g}(\integ_\ell)$; still,
\cite{fw} entertains the hope that this distinction is immaterial to
the problem.)  Friedman and Washington observe that, since geometric
equidistribution results seem within reach, their conjecture may well
be tractable.

Taking advantage of recent progress in equidistribution, we
prove a statement in the spirit of \cite{fw}.  A special, yet typical,
case of our main result  says the following.

Let $\calc \ra \calm/\ff_q$ be a relative smooth, proper curve of
genus $g$ over a smooth, irreducible variety, and suppose that $\calc
\ra \calm$ has full $\ell$ monodromy.  Let $\st{\ff_{q^{e_n}}}$ be a
tower of finite extensions of $\ff_q$ which is 
cofinal in the collection of all finite extensions of $\ff_q$.  For a
curve $C$, let $\jac(C)[\ell](k)$ denote the $k$-rational $\ell$-torsion subgroup of its Jacobian.
We give an explicit formula for a number $\alpha(g,r)$ so that:
%\begin{theorem}
\[
\lim_{n\ra\infty}\frac{\abs{\st{ x \in \calm(\ff_{q^{e_n}}) :  
\jac(\calc_x)[\ell](\ff_{q^{e_n}}) \iso
(\integ/\ell)^r}}}{\abs{\calm(\ff_{q^{e_n}})}} = \alpha(g,r).
\]
%\end{theorem}

Moreover, $\lim_{g \ra \infty} \alpha(g,r)$ exists.  In this way, we
can formulate a version of our result which allows the genus of the
curves in question to change, too. The term $\alpha(g,r)$ should be
thought of as a sort of symplectic analogue of $\abs{\aut((\integ/\ell)^r)}\inv$.

This result lets us essentially prove the original conjecture
of Friedman and Washington.  (We will see in Section \ref{subsecfw} that the heuristic
they introduce, that statistics of $\gl_{2g}$ should track those of
$\sp_{2g}$, is a reasonable approximation but not literally true.)
Moreover, we significantly strengthen (Section \ref{subseche})
results of Cardon and Murty \cite{cardonmurty} on divisibility of
class groups of quadratic function fields.

The family of curves $\calc \ra \calm$ is a concrete device for
enumerating function fields.  The ``full $\ell$-monodromy''
constraint ensures that, as far as $\ell$-Sylow subgroups of class
groups are concerned, the family behaves like a general one.

On one hand, the familiar moduli spaces $\mgn$ \cite{dm} of 
proper smooth curves of genus $g$ equipped with principal Jacobi level
$N$ structure
 have
this property.  In fact, so do most versal families of curves
\cite{ekedahl}.  In this sense, our main theorem describes a typical
collection of function fields of genus $g$.

On the other hand, it's not hard to write down families of curves
which {\em do not} have this property.  Generally speaking, extra
algebraic cycles on a family of curves force the $\ell$-adic monodromy
to lie in a proper subgroup of $\sp_{2g}$.  As a specific caution, we
mention that if $d > 2$ then the family of curves $C_{d,f}: y^d =
f(x)$ does not have full monodromy.  (In fact, for such curves,
$\integ[\zeta_d]  \inject \End_{\bar\ff_q}\jac(C_{d,f})$.  One
suspects \cite{zarendjac} that the $\ell$-adic monodromy is a unitary
group associated to $\rat(\zeta_d)$; at the very least, the monodromy
group is contained in such a group.)

The first section of this paper collects results of Katz about
equidistribution of Frobenius elements on $\ell$-adic sheaves.  The
second section investigates the combinatorics of $\sp_{2g}(\ff_\ell)$
and related groups.  The next  section combines these results to make
precise statements about the distribution of class groups of function
fields.   The paper concludes with a series of applications of these
results, culminating in a proof of a modified Friedman-Washington
conjecture.

I thank R. Pries for helpful comments on this paper.

\section{$\ell$-adic monodromy}
\setcounter{equation}{0}

As noted above, Friedman and Washington foresaw that good
equidistribution theorems would allow one to prove Cohen-Lenstra type
results for function fields.  Here, we recall the precise statements
we need.  Our discussion follows section one of \cite{achhol}, which
itself is a recapitulation of parts of chapter nine of
\cite{katzsarnak}.

Fix an odd prime $\ell$.  Let $\calo_\lambda$ be the ring of integers
in some finite extension of $\integ_\ell$, and let $\Lambda =
\calo_\lambda/\lambda^n$ for some $n$.  Let $V = V_\Lambda$ be a free,
rank $2g$ $\Lambda$-module equipped with a symplectic form
$\ip{\cdot,\cdot}$.  The group of symplectic similitudes of
$(V,\ip{\cdot,\cdot})$ is
\[
\gsp(V,\ip{\cdot,\cdot}) = \st{
A \in \gl(V) |  \exists \mult(A)\in \Lambda\units  : \forall v,w\in V,
\ip{Av,Aw} = \mult(A)\ip{v,w}} \iso \gsp_{2g}(\Lambda).
\]

The ``multiplicator'' $\mult$ is a character of 
$\gsp(V,\ip{\cdot,\cdot})$, and its
kernel is the usual symplectic group $\sp_{2g}(\Lambda)$.  For
$\xi\in \Lambda\units$, let $\gsp_{2g}^\xi(\Lambda) =
\mult\inv(\xi)$ be the 
set of symplectic similitudes with multiplier $\xi$; each
$\gsp_{2g}^\xi$ is a torsor over $\sp_{2g}$.  For $W\subset
\gsp_{2g}(\Lambda)$, let $W^\xi = W\cap \gsp_{2g}^\xi(\Lambda)$.

Let $T \ra \spec \integ[1/\ell]$ be a connected normal scheme of
finite type, often $\spec\integ[1/\ell]$ itself. Let $\calm \ra T$ be
a scheme with smooth geometrically irreducible fibers; let
$\eta_\calm$ be a generic point of $\calm$.  A local system $\calf$
of symplectic $\Lambda$-modules of rank $2g$ is equivalent to a continuous
representation $\rho_\calf:\pi_1(\calm,\bar\eta_\calm) \ra
\aut(\calf_{\bar\eta_\calm}) \iso \gsp_{2g}(\Lambda)$.  We call the image of
this representation the (arithmetic) monodromy group of $\calf$.

Let $k$ be a finite field and $t\in T(k)$.  Then $\calm_t$ is a
$k$-scheme, and we distinguish the geometric fundamental group
$\pi_1^\geom(\calm_t) = \pi_1(\calm_t \cross \bar k)\subseteq \pi_1(\calm_t)$.   The Galois
group of $k$ is (canonically isomorphic to) the quotient
$\pi_1(\calm_t)/\pi_1^\geom(\calm_t)$. 

We will require our sheaves to have uniform geometric monodromy group
$G_\geom\subseteq \sp_{2g}(\Lambda)$, in the sense that for every
finite field $k$ and $t\in T(k)$, the image of
$\rho_\calf(\pi_1^\geom(\calm_t))$ is $G_\geom$.  With this assumption
the monodromy group $G$, the full image of $\rho_\calf$, is contained
in $\Lambda\cdot G_\geom \subseteq \gsp_{2g}(\Lambda)$.  We let
$\xi(k)$ denote the image of $\frob_{k}$, the canonical generator of
$\pi_1(\spec k)$, in $G/G_\geom \subseteq \Lambda\units$.

If $k$ is a finite field, then to a $k$-point $x\in \calm(k)$
one may associate its (conjugacy class of)
Frobenius $\frob_{x/k}$ in $\pi_1(\calm)$.  Via $\rho_\calf$, the
Frobenius at $x$ acts on $\calf_{\bar\eta_\calm}$.  Katz shows that these
Frobenius elements are equidistributed in the monodromy group.

\begin{theorem}[Katz]\label{katz} Suppose $\calf$ has uniform geometric
monodromy group $G_\geom$ and arithmetic monodromy group $G$.    Let
$W\subset G$ be stable under $G$-conjugation.  There are effective
constants $\delta(\calm,\calf)$ and $A(\calm/T)$ so that, if $k$ is a
finite field with $\abs k >A(\calm/T)$ and $t:\spec k \ra T$ is an
inclusion, then
\[
\abs{
\frac{\abs{\st{ x\in \calm_t(k): \rho_\calf(\frob_{x,k})\in
      W}}}{\abs{\calm_t(k)}}-
\frac{\abs{W^{\xi(k)}}}{\abs{G^{\xi(k)}}}
}
< \epsilon(\calm,\calf,k) \defeq \frac{\delta(\calm,\calf)}{\sqrt{\abs k}}.
\]
\end{theorem}

\pf This is simply \cite[9.7.13]{katzsarnak}; see also
\cite[4.1]{chav}. While the result holds for any sheaf with finite
monodromy group, we will (almost; see Section \ref{subsececqn1}
below) always work with subgroups of $\gsp_{2g}(\Lambda)$.
\qed

Already, this deep theorem yields a method for computing the 
proportion of curves in a family for which the $\ell$-Sylow subgroup 
of the class group is isomorphic to a given group.  Indeed, let $H$ be 
any finite abelian group annihilated by $\ell^e$, and let
$\pi:\calc\ra\calm/T$ be a smooth, irreducible proper relative curve of
genus $g\ge 1$.   For any finite field $k$ and $t\in T(k)$ we define
\begin{align}\label{defbeta}
\beta(\calc\ra\calm, t, \ell^e,H) & =  
\frac{\abs{\st{x \in \calm_t(k): \jac(\calc_x)[\ell^e](k) \iso H}}}
{\abs{\calm(k)}}.
\end{align}
We will often assume the $k$-point $t:\spec k \ra T$ is fixed, and
simply write $\beta(\calc\ra\calm, k, \ell^e, H)$.

There is a sheaf $\calf = \calf_{\calc,\ell^e}$ of abelian groups on
$\calm$ whose fiber at a geometric point $\bar x\in \calm$ is the
$\ell^e$-torsion of the Jacobian $\jac(\calc_{\bar x})[\ell^e]$; it
may be alternatively defined by $\calf = R^1\pi_!(\integ/\ell^e)$.
Suppose that this family has uniform geometric monodromy group
$G_\geom \subseteq \sp_{2g}(\integ/\ell^e)$ and arithmetic monodromy
group $G\subseteq \gsp_{2g}(\integ/\ell^e)$.  Then Theorem \ref{katz} implies that, for any
sufficiently large finite field $k$ and fixed, suppressed $k$-point of
$T$,  we have
\begin{align}\label{defepsilon}
\abs{\beta(\calc\ra\calm,k,\ell^e,H)-\frac{\abs{\st{ x\in G^{\xi(k)} : \ker(x-\id) \iso
H}}}{\abs{G^{\xi(k)}}}} &< \epsilon(\calm,\calf,k).
\end{align}
We will denote the right-hand term of this inequality by
$\epsilon_{\calc\ra\calm}(\ell^e,k)$.  Note that for a fixed family of
curves $\calc\ra \calm$, as $\abs k \ra \infty$ we have
$\epsilon_{\calc\ra\calm}(\ell^e,k) \ra 0$.

In the special case where $\ell H = 0$, the mod-$\ell$ monodromy group
is the full symplectic group, and $\st{k_n}$ is a collection of finite
fields equipped with maps to $T$; $\lim_{n\ra\infty}\abs{k_n} =
\infty$; 
and, for $n\gg0$, $\abs{k_n}\equiv 1 \bmod \ell$;  we have 
\begin{align}\label{eqfullmono} 
\lim_{n\ra\infty} \beta(\calc\ra\calm, k_n, \ell, H) &= 
\frac{\abs{\st{ x \in \sp_{2g}(\ff_\ell) : \ker(x - \id) \iso 
      H}}}{\abs{\sp_{2g}(\ff_\ell)}}. 
\end{align} 
In the next section we explain how to compute the right-hand side of 
\eqref{eqfullmono}. 

We collect the diverse notation and assumptions of this section in the
following:

\begin{situation}\label{mysit} We suppose that $\ell$ is a fixed, odd prime;
$T$ is a
connected $\integ[1/\ell]$-scheme of finite type; $\calm \ra \calt$ 
is a smooth scheme with geometrically irreducible fibers; $\calc \ra
\calm$ is a proper, smooth relative curve of genus $g$;
$\calf_{\ell^e}$ is the sheaf of $\ell^e$-torsion on the Jacobian of
$\calc$; and $\calf$ has uniform geometric monodromy group $G_\geom$ and
arithmetic monodromy group $G$.  We will say that $\calc \ra \calm$
has full $\ell^e$-monodromy if $G_\geom = \sp_{2g}(\integ/\ell^e)$.
Additionally, $\st{k_n}$ is a collection of finite fields, each
equipped with an inclusion $t_n:k_n \ra T$, such that
$\lim_{n\ra\infty}\abs{k_n} = \infty$, and we will often write
$\calm(k)$ for $\calm_{t_n}(k)$.
\end{situation}

While any sequence of finite fields is allowed, psychologically it
seems to be easiest to think of either $\st{\ff_{q^{e_n}}}$, a tower of
extensions of a fixed finite field $\ff_q$, or $\st{\ff_{p_n}}$, a
collection of finite fields of ever-larger prime order.

\section{Matrices with given fixed space}\label{secsp}

\def\sl{{\rm SL}}
\def\so{{\rm SO}}

Katz's work on mod-$\ell$ monodromy reduces the calculation of $\beta$
to a calculation in a symplectic group $\sp_{2g}(\ff_\ell)$.  In
this section we go to some length to calculate precisely the
proportion of elements with a certain behavior.  We note that \cite{fw}
avoids these difficulties in two ways.  First, the authors compute in
$\gl_n$, where the relevant Lie theory is more transparent,
rather than in $\sp_{2g}$.  Second, they compute with
elements of $\mat_n$, with the hope that if $x$ is equidistributed in
$\gl_n$, then $x-\id$ is equidistributed in $\mat_n$.  Unfortunately,
these choices mean that the conjectural description of class group
frequencies in \cite{fw} differs slightly from the actual frequencies;
we take up this point in more detail in Section \ref{subsecfw}.

\subsection{Symplectic matrices over finite fields}\label{subsecsymp}

Fix a finite field $\ff$ with $\ell$ elements, where $\ell$ is a power
of an odd prime.  (In our applications $\ell$ will itself be prime,
but this assumption is not necessary for the present computation.)

Our main goal is a formula for 
\begin{eqnarray}\label{defalpha}
\alpha(g,r)&\defeq& \frac{\abs{\st{x\in \sp_{2g}(\ff) : \ker(x-\id)\iso \ff^r}}}{\abs{\sp_{2g}(\ff)}}.
\end{eqnarray}
The method presented here should work for any of the classical
families of finite groups of Lie type.  Still, since it is the symplectic group which arises
most naturally in questions about the typical function field, we have
chosen to focus our efforts on groups of type $C$.  The reader will
notice our heavy reliance on the paper \cite{springersteinberg} of
Springer and Steinberg.

Consistent with the notation introduced in the previous section, we
view $\sp_{2g}(\ff)$ as the group of automorphisms of a
$2g$-dimensional $\ff$-vector space $V_g$ equipped with a symplectic
form $\ang{\cdot,\cdot}_g$.   Unless otherwise noted, an $r$-subspace of
$V_g$ means any subspace $W\subset V_g$ for which
$(W,\ang{\cdot,\cdot}_g\rest W)\iso (V_r, \ang{\cdot,\cdot}_r)$.  
We will need to isolate the subspace of $V_g$ on which a given
element $x\in \sp(V_g)$ acts unipotently.

\begin{lemma}\label{eigenspace} Suppose $x\in \sp(V)\iso \sp_{2g}(\ff)$.  Then there are subspaces $E_1(x)$
and $E_1(x)^\perp$ such that $V\iso E_1(X) \oplus E_1(x)^\perp$;
$x\rest{E_1(x)}$ is unipotent; and $x-\id$ is invertible on
$E_1(x)^\perp$.
\end{lemma}

%% \begin{lemma} Suppose $x\in \sp(V)\iso \sp_{2g}(\ff)$.  Then there are subspaces $E_1(x)$
%% and $E_1(x)^\perp$ such that $V\iso E_1(X) \oplus E_1(x)^\perp$;
%% $x\rest{E_1(x)}$ is unipotent; and $x-\id$ is invertible on
%% $E_1(x)^\perp$.
%% \end{lemma}

\pf We assume that $x-\id$ is not invertible, as otherwise the statement is trivial.
Write $V_{\bar\ff}$ as the direct sum of generalized eigenspaces
for $x$, $V_{\bar\ff} = \oplus V_{\bar\ff}(\lambda)$ where
$V_{\bar\ff}(\lambda)$ is the kernel of $(x-\lambda\id)^{2g}$ on
$V_{\bar\ff}$.

Suppose that $\lambda\mu\not = 1$. We will prove, by induction on
$m+n$, that 
\[
\ip{\ker(x-\lambda\id)^m\rest{V_{\bar\ff}},
  \ker(x-\mu\id)^n\rest{V_{\bar\ff}}} = 0.
\]
For the base case $m=n=1$,
suppose $xu = \lambda u$ and $x v = \mu v$.  Because $x$ preserves the
symplectic form, we have  $\ip{u,v} =
\ip{xu,xv} = \ip{\lambda u, \mu v} = \lambda\mu\ip{u,v}$.  Since
$\lambda\mu\not = 1$, this forces $\ip{u,v} = 0$.

We now treat the inductive step.  Suppose that $u\in
\ker(x-\lambda\id)^m \rest{V_{\bar\ff}}$ and $v\in \ker(x-\mu\id)^n
\rest{V_{\bar\ff}}$.  Without loss of generality, assume that $m\ge n
\ge 1$.  Then $xu = u' + \lambda u$, where $u' \in
\ker(x-\lambda\id)^{m-1}\rest{V_{\bar\ff}}$, and $xv = v'+\mu v$ for
some $v'\in \ker(x-\mu\id)^{n-1}\rest{V_{\bar\ff}}$.  (If $n=1$, this
simply means that $v'=0$.)  We then have $\ip{u,v} =
\ip{xu,xv} = \ip{u'+\lambda u, v' + \mu v} = \ip{u',v'}+\lambda
\ip{u,v'} + \mu\ip{u',v} + \lambda\mu\ip{u,v}$.  By the inductive
hypothesis, the first three terms in the last expression vanish.  This
leaves us with $\ip{u,v} = \lambda\mu\ip{u,v}$; again, $\ip{u,v} = 0$.

This shows that, if $\lambda\mu\not =1$, then
$\ip{V_{\bar\ff}(\lambda), V_{\bar\ff}(\mu)} = 0$.  Since the pairing
$\ip{\cdot,\cdot}$ is nondegenerate on $V$, we conclude that the pairing
$\ip{\cdot,\cdot}: V_{\bar\ff}(\lambda) \cross
V_{\bar\ff}(\lambda\inv) \ra \bar\ff$ is nondegenerate.  In
particular, $V_{\bar\ff}(1)$ is self-dual under $\ip{\cdot,\cdot}$.

Now, the generalized eigenspace associated to $1$ is defined over
$\ff$; therefore, its orthogonal complement is, too.    Returning to
the $\ff$-vector space $V$, we find that $E_1(x) \defeq
\ker(x-\id)^{2g}\subset V$ is a symplectic subspace of $V$.
Therefore, there exists a canonical decomposition $V = E_1(x)
\oplus E_1(x)^\perp$, where $1$ is not an eigenvalue of the action of
$x$ on $E_1(x)^\perp$.\qed

We define the following quantities associated to $\sp_{2g}(\ff)$.  Let
$\nu(g)$ be the number of elements in $\sp_{2g}(\ff)$; let $U(g)$ be
the number of unipotent elements in $\sp_{2g}(\ff)$; and let $S(g,r)$
be the number of $r$-subspaces of $V_g$.  Let $\Phi(g)$ be the number
of elements $x\in\sp_{2g}$ for which $x-\id$ is invertible, and let
$\phi(g) = \Phi(g)/\nu(g)$ be the proportion of symplectic matrices
with this property.  For convenience, we define $\Phi(0) = 1$.

Recall that $\alpha(g,r)$ is the proportion of elements $x\in
\sp_{2g}(\ff)$ for which $\ker(x-\id) \iso \ff^r$.  Let $U(g,r)$ be
the number of {\em unipotent} elements $u$ of $\sp_{2g}(\ff)$ for
which $\ker(u-\id) \iso \ff^r$.  These quantities enjoy the following
relations.

\begin{lemma}\label{calcphi} With all notation as above, let $\lambda(j) =
\ell^{2j-1}(\ell^{2j}-1)$.  Then $\nu(g) =
\prod_{j=1}^g \lambda(j)$; $S(g,r) = \nu(g)/(\nu(r)\nu(g-r))$; $U(j) =
\ell^{2j^2}$; 
\begin{align}
\alpha(g,r) &= \oneover{\nu(g)}\sum_{j=1}^g S(g,j)
U(j,r)\Phi(g-j)\label{eqcalcalpha};\\
\intertext{and}
\Phi(g)& = \nu(g) - \sum_{j=1}^g S(g,j)U(j) \Phi(g-j).\label{eqcalcphi}
\end{align}
\end{lemma}

\pf The calculation of $\nu$ and $S$ is standard geometric algebra
\cite[III.6]{artingeomalg}.  One proves that the symplectic group acts simply
transitively on symplectic bases for $V_g$, and that $\lambda(g)$
counts the number of symplectic pairs in $V_g$.  A theorem of
Steinberg (\cite[8.14]{humphreysconjclass} or \cite{springersteinberg}) says that the number of unipotent elements in a finite group
$G$ of Lie type is $\ell^{\dim G - \rank G}$.  Therefore $U(g)$, the
number of unipotent elements in $\sp_{2g}(\ff)$, is $\ell^{2g^2}$.

By Lemma \ref{eigenspace}, any $x\in
\sp_{2g}(\ff)$ determines a decomposition $V = E_1(x) \oplus
E_1(x)^\perp$, where $x$ acts unipotently on $E_1(x)$ and $(x-\id)$ is
invertible on $E_1(x)^\perp$.
Therefore, any element of the
symplectic group determines, and is determined by, the data of a
subspace $W\subset V$; a unipotent element $u\in \sp(W)$; and an
element $y\in \sp(W^\perp)$ for which $(y-\id)$ is invertible. If $x$
corresponds in this way to the triple $(W,u,y)$, then $\ker(x-\id) =
\ker(u-\id)\rest W$.

Equations \eqref{eqcalcalpha} and \eqref{eqcalcphi} follow swiftly.
The right-hand side of \eqref{eqcalcalpha} enumerates all choices of
data $(W,u,y)$ where $W$ is a $j$-subspace of $V$, $u$ is a unipotent
element of $\sp(W)$ with $\ker(u-\id) \iso \ff^r$, and $y\in
\sp(W^\perp)$ with $y-\id$ invertible, all normalized by the size of
the symplectic group.

  To calculate $\Phi(g)$ and thus derive \eqref{eqcalcphi}, we simply subtract from
$\nu(g)$ the number of symplectic elements with nontrivial unipotent
part.  We enumerate triples $(W,u,y)$ as before, where $W$ is a
positive-dimensional subspace of $V$.  If $W \iso V_j$, then $U(j)$
counts the number of choices for $u$, while $\Phi(g-j)$ is, by
definition, the number of choices for $y$.\qed

Equation \eqref{eqcalcalpha}, combined with Proposition \ref{calcu}
below, allows the explicit computation of $\alpha(g,r)$ in any
particular case.  The results of this calculation for $g \le 3$ are
shown in Table \ref{bigtable}.

\begin{prop}\label{calcu} The number of unipotent elements $u$ in $\sp_{2g}(\ff)$ such
that $\ker(u-\id)\iso \ff^r$ is
\begin{align}
U(g,r) &= \nu(g) \sum_{{\bf d}: 0 < d_1 \le d_2 \le \cdots \le d_r}
\left(
\ell^{\half 1 \left( \sum_i s_i^2 - \sum_i r_i^2 + \sum_{i\text{ even}}
r_i\right)}
\cdot 
\prod_{i\text{ odd}} \nu(r_i/2)
\cdot
\prod_{i\text{ even}}\nu_\orth(r_i)
\right)\inv\label{eqbigmess}
\end{align}
where the sum is over all partitions ${\bf d}$ of $\dim V_g$ into $r$
parts such that odd parts occur with even multiplicity; $r_i = \abs{\st{
  j : d_j =i }}$; $s_i= \sum_{j\ge i} r_i$; and
\[
\nu_\orth(n) = \begin{cases}
\ell^{m^2}\prod_{i=1}^m (\ell^{2i}-1) & n = 2m+1 \\
\ell^{m^2-2m}\prod_{i=1}^{m}(\ell^{2i}-1) & n = 2m
	       \end{cases}
.
\]
\end{prop}

\pf Let $G = \sp_{2g}$.  We start by identifying the relevant
$G(\bar\ff)$-conjugacy classes of unipotent elements.
As in, say, \cite[6.20]{humphreysconjclass}, 
let
$\calu = \calu(G)$ be the unipotent variety of $G$; it parametrizes
all unipotent elements of $G$.  Similarly, let $\caln$ be the
nilpotent variety of $\lie g$, the Lie algebra of $G$.  The Cayley
transform is a $G$-equivariant isomorphism
\begin{equation*}
\xymatrix{
\calu  \ar[r]&  \caln\\
x  
\ar@{|->}[r] & (1-x)(1+x)\inv.
}
\end{equation*}

Thus, it suffices to count those $y\in \caln(\ff)$ with
 nullspace of rank $r$.

Happily, enumeration of nilpotent elements is a classical result.
Moreover, the description makes it easy to pick out those with the
appropriate rank.
To give a nilpotent orbit in $\lie{s}\lie{l}_n$ is to describe its Jordan
normal form; a similar classification exists for arbitrary Lie groups.  We have
the classical bijection
(\cite[5.1.1]{collingwoodmcgovern}, \cite[7.11]{humphreysconjclass})
between nilpotent orbits of $\lie g$ and the partitions of $2g$ for
which odd parts occur with even multiplicity.
The dimension of the nullspace of an element in a nilpotent
orbit corresponding to a given partition is the number of elements in
that partition.  Therefore, the desired (geometric) nilpotent orbits are
represented by suitable partitions with exactly $r$ pieces.   Each of
these conjugacy classes has a representative in $G(\ff)$, and 
the summation in equation \eqref{eqbigmess} thus ranges over all
$G(\bar\ff)$-conjugacy classes of unipotent elements $x$ in $G(\ff)$
for which $\ker(x-\id) \iso \ff^r$.

We now explain how these $G(\bar\ff)$ conjugacy classes behave
over $G(\ff)$, and compute the isomorphism class of the
centralizer (still in $G(\ff)$) of an element of such a
conjugacy class.  By doing so, we are able to compute the size of the
relevant conjugacy class.

We proceed as in \cite[IV.2]{springersteinberg}.  Fix a geometric
conjugacy class corresponding to a partition ${\bf d}$ of $g$, and let
$I = I({\bf d})= \st{ i : i\text{ even and }r_i >0 }$.  Jordan factors
corresponding to even members $d_i$ split into two conjugacy classes
over $\ff$.  Therefore, to give a $G(\ff)$-conjugacy class inside the
$G(\bar\ff)$ conjugacy class ${\bf d}$ is to give a map of sets $c: I
\ra \st{-1,+1}$.

Let $u$ be a representative for the $G(\ff)$-conjugacy class
corresponding to ${\bf d}$ and a choice of assignments $c$.   A
theorem of Springer and Steinberg \cite[IV.2.26-8]{springersteinberg} computes the isomorphism class of
the centralizer $Z = Z_G(u)$ in $G$.  It is the semidirect product of
a unipotent radical, $R$, and the centralizer $C$ of a certain
torus associated to $u$.  (Note that \cite{springersteinberg} computes
the connected component of the centralizer, and then later accounts
for multiple components.)  The dimension of the Lie algebra
of $R$ is $\half 1 \left( \sum_i s_i^2 - \sum_i r_i^2 + \sum_{i\text{ even}}
r_i\right)$.  The reductive group $C$ is isomorphic to
\[
\prod_{i\text{ odd}} \sp_{r_i}(\ff) \cdot \prod_{i\text{ even}}
\O^{c(i)}_{r_i}(\ff).
\]
Here, if $r_i$ is even then $\O^{+1}_{r_i}(\ff)$ denotes the rank $r_i$
orthogonal group of Witt defect $0$ over $\ff$, while
$\O^{-1}_{r_i}(\ff)$ is the orthogonal group of Witt defect $1$.  For
odd $r_i$, $\O^{\pm 1}_{r_i}(\ff)$ is the (unique) orthogonal group of
rank $r_i$.

Therefore, the size of the set of elements in $G(\ff)$ which belong to
the $G(\bar\ff)$-conjugacy class represented by ${\bf d}$ is
\begin{align*}
%% \frac{\abs{\sp_{2g}(\ff)}}{\abs R \prod_{i\text{ odd}}
%%   \abs{\sp_{r_i}(\ff)}} \cdot \left(
%% \sum_{c: I \ra \st{-1,+1}} \prod_{i\in I}
%%   \abs{\O_{r_i}^{c(i)}(\ff)}\inv \right) &=
\sum_{c:I \ra \st{-1,+1}} 
\frac{\abs{\sp_{2g}(\ff)}}{\abs R \prod_{i\text{ odd}}
  \abs{\sp_{r_i}(\ff)}}
\prod_{i\in I}\abs{\O_{r_i}^{c(i)}(\ff)}\inv
&=
\frac{\nu(g)}{
\abs R \cdot
%\ell^{\half 1 \left( \sum_i s_i^2 - \sum_i r_i^2 + \sum_{i\text{ even}}
%r_i\right)}
 \prod_{i\text{ odd}}\nu(r_i/2)} \cdot \prod_{i\in I}
  \left(
\abs{\O_{r_i}^{(-1)}(\ff)}\inv + \abs{\O_{r_i}^{(+1)}(\ff)}\inv\right)
  \\
&=
\frac{\nu(g)}{
\ell^{\half 1 \left( \sum_i s_i^2 - \sum_i r_i^2 + \sum_{i\text{ even}}
r_i\right)} \prod_{i\text{ odd}}\nu(r_i/2)} \cdot \prod_{i\in I}
  \nu_\orth(r_i),
\end{align*}
where
\[
\nu_\orth(n) = \begin{cases}
\ell^{m^2}\prod_{i=1}^m (\ell^{2i}-1) & n = 2m+1 \\
\ell^{m^2-2m}\prod_{i=1}^{m}(\ell^{2i}-1) & n = 2m
	       \end{cases}
.
\]
Note that $\nu_\orth(2m+1)$ is simply the number of elements in
$\so_{2m+1}(\ff)$, while $\nu_\orth(2m+1)$ is the harmonic mean of
$\abs{\so^{(-1)}_{2m}(\ff)}$ and $\abs{\so^{(+1)}_{2m}(\ff)}$.  By
summing over suitable geometric conjugacy classes ${\bf d}$ we obtain
equation \eqref{eqbigmess}.\qed

\begin{lemma}\label{alphainf} The limits
\[
\phi(\infty)\defeq \lim_{g\ra\infty} \phi(g)\text{ and
}\alpha(\infty,r) \defeq \lim_{g\ra\infty}\alpha(g,r)
\]
exist.
\end{lemma}

\pf By Lemma \ref{calcphi},
\begin{eqnarray*}
\phi(g) & = & \oneover{\nu(g)}(\nu(g) - \sum_{j=1}^g S(g,j)U(j) \Phi(g-j))
\\
& = & 1 - \sum_{j=1}^g \frac{U(j)\Phi(g-j)}{\nu(j)\nu(g-j)}.
\end{eqnarray*}
Now, $\Phi(g-j)$ is necessarily less than $\nu(g-j)$, while
$U(j)/\nu(j) < \ell^{-j}$.  Therefore, $\lim_{g\ra\infty}\phi(g)$
exists.  Similarly, consider
\[
\alpha(g,r) = \sum_{j=1}^g \frac{U(j,r)\Phi(g-j)}{\nu(j)\nu(g-j)}.
\]
Again, $U(j,r)/\nu(j) \le U(j) / \nu(j) < \ell^{-j}$, so that $\lim_{g\ra\infty}\alpha(g,r)$ converges.\qed

%% \subsection{Finite abelian $\ell$-groups}
%% Given an arbitrary finite
%% abelian $\ell$-group $H$, let
%% \[
%% \alpha(g,H,\ell^e) = \frac{\abs{\st{x\in \sp_{2g}(\integ/\ell^e)
%% :\ker(x-\id) \iso H \bmod \ell^e}}}{\abs{\sp_{2g}(\integ/\ell^e)}}.
%% \]
%% While Section \ref{subsecsymp} explains how to compute this quantity
%% in the special case where $e=1$, we do not know a general answer.
%% Still, in any particular case it is not hard to calculate this value,
%% and we can give the following brief indications.

%% First, $\alpha(g,H,\ell^n)$ stabilizes as $n$
%% gets large.  Indeed, if $\ell^e H = 0$ then $\alpha(g,H,\ell^e) = \alpha(g,H,\ell^{e+m})$,
%% since the projection $\sp_{2g}(\integ/\ell^{e+n}) \ra
%% \sp_{2g}(\integ/\ell^e)$ is a surjection of groups.

%% Second, the limit $\lim_{g\ra\infty}\alpha(g,H,\ell^e)$ exists.  If we
%% define quantities such as $U(j,r,\ell^e)$ in the expected way, they
%% are sufficiently close to $U(j,r,\ell) = U(j,r)$ that the proof of
%% Lemma \ref{alphainf} goes through.

\subsection{Unitary groups}\label{subsecunitary}
The methods of Section \ref{subsecsymp} work for any family of
classical Lie groups.  Since unitary groups also come up in certain natural
applications (see Section \ref{extrigonal}), we briefly indicate how the
argument works for $\u_n$.  Because $\u_n$ is a twist of $\gl_n$, the
details are actually somewhat simpler.  We preserve all notation from Section
\ref{subsecsymp}, using the subscript $\u$ to denote the appropriate
group.

So, let $\u_n$ denote the unitary group in $n$ variables over $\ff$.
Implicit in this definition is a nontrivial involution $\sigma$ of
$\ff$; let $m$ be $\sqrt{\ell}$, the size of the fixed field of $\sigma$.
The number of elements in $\u_n$ is
\[
\nu_\u(n) = m^{\half 1(n^2-n)} \prod_{i=1}^n(m^i - (-1)^i);
\]
the number of unipotent elements is $U_\u(n) = m^{n^2-n}$; and
$S_\u(n,r) = \nu_\u(n)/(\nu_\u(r)\nu_\u(n-r))$.  Moreover, the
number of unitary matrices for which $1$ is not an eigenvalue is
\[
\Phi_\u(n) = \nu_\u(n) - \sum_{j=1}^n S_\u(n,j)U_\u(j)\Phi_\u(n-j).
\]
Since the unitary group is a form of the general linear group,
unipotent classes are parametrized by (unrestricted) partitions of
$n$.  Moreover, $\u_n(\bar\ff)$ and $\u_n(\ff)$ conjugacy coincide.  The
centralizer of a unipotent element corresponding to the partition $0 <
d_1 \le \cdots \le d_r$ of $n$ is connected, and has
\[
m^{\sum s_i^2 - \sum r_i^2} \prod \nu_\u(d_i)
\]
elements.  All other results of Section \ref{subsecsymp}, including
the existence of $\alpha_\u(\infty, r)$, carry over.

\section{Class groups of families of curves}

Let $H$ be a finite abelian group with $\ell^e H =0$.  We would like
to compute the chance that $H$ is the $\ell$-part of the class group
of a function field.  As ``sample space'' of function fields we choose
the fibers of
any relative curve $\calc \ra \calm/T$  as in Situation \ref{mysit}
with full $\ell^e$-monodromy.

 In practice, general families
of curves tend to have full $\ell^e$-monodromy; see, for instance, the
introduction to \cite{ekedahl}.  As a concrete example, fix a natural
number $N\ge 3$ relatively prime to $p$ and consider $\cgn\ra \mgn$,
the universal curve of genus $g$ with principal Jacobi structure of
level $N$.  By \cite[5.15-5.16]{dm}, this family of curves has
full $\ell^e$ monodromy.  Indeed, any versal family of curves has full
monodromy at most primes \cite[2.2]{achhol}.  We also expect
(see Section \ref{subseche}) that a general family of
hyperelliptic curves has full $\ell^e$-monodromy.

The equidistribution results in the first section let us
detect the occurrence of $H$ in class groups of function fields.
Recall (Equation \eqref{defbeta}) that this is measured by
$\beta(\calc\ra\calm,k, \ell^e,H)$.

\begin{theorem}\label{maintheorem} Let $H$ be a finite abelian $\ell$-group.
As in Situation \ref{mysit}, let $\calc \ra \calm/T$ be a relative
curve with full $\ell^e$-monodromy and let $\st{k_n}$ be a collection
of finite fields.  Suppose that for $n\gg 0$, $\xi(k_n) \equiv \xi \bmod
\ell^e$.  There exists an effective constant $\delta(\calc \ra \calm)$
so that, for $n$ sufficiently large,
\begin{align*}
\abs{\beta(\calc\ra\calm, k_n, \ell^e, H) - \alpha^\xi(g,H,\ell^e)} &< \epsilon_{\calc \ra \calm}(\ell^e,k_n) \defeq
  \frac{\delta(\calc\ra\calm)}{\sqrt{\abs{k_n}}},\\
\intertext{and thus }
\lim_{n\ra\infty}\beta(\calc\ra\calm,k_n,\ell^e,H)& = \alpha^\xi(g,H,\ell^e),
\end{align*}
where
\[
\alpha^\xi(g,H,\ell^e) = \frac{\abs{\st{x \in
\gsp^\xi_{2g}(\integ/\ell^e) : \ker( x-\id) \iso H
}}}{\abs{\sp_{2g}(\integ/\ell^e)}}.
\]
In the special case where $e = 1$ and $\xi = 1$, this term is computed
by Lemma \ref{calcphi} and Proposition \ref{calcu}.
\end{theorem}

\pf As above consider the lisse sheaf $\calf = \calf_{\calc,\ell^e}$
on $\calm$ which associates, to each geometric point $\bar x$, the 
$\ell^e$-torsion of the Jacobian of $\calc_{\bar x}$.  Let  $x\in \calm(k_n)$
be any point.  The divisor class group of $\calc_x$ is
$\jac(\calc_x)(k_n)$, the $k_n$-rational points of the Picard variety.
By definition, the $\ell^e$-torsion of this group is the subgroup of $\calf_x$ fixed by $\frob_{x/k_n}$.  Thus, in the
notation of the first section, $\jac(\calc_x)[\ell^e](k_n) \iso
\ker (\rho_\calf(\frob_{x/k_n})-\id)$, and
\[
\beta(\calc\ra\calm, k_n, \ell^e, H)  = \frac{\abs{\st{ x \in \calm(k_n) : \ker
(\rho_\calf(\frob_{x/k_n})-\id) \iso H}}}{\abs{\calm(k_n)}}.
\]
In general, Theorem \ref{katz} finishes the proof.  For the special case where
$H$ is an elementary abelian $\ell$-group and $\abs{k_n}\equiv 1 \bmod \ell$,
Section \ref{secsp}
provides an algorithm for computing the appropriate quantity.\qed

%% Note that if $\ell^{e-1}H$, then Theorem \ref{maintheorem} states that
%% $\alpha(g,H,\ell^e)$ is the proportion of curves for which the
%% $\ell$-Sylow subgroup of the class group is isomorphic to $H$.

In some applications, it is useful to be able to consider a family of
curves with unbounded genus.   To employ our methods, we
need the size of the field of constants to grow more swiftly than the
error terms $\epsilon$ of \eqref{defepsilon}.

\begin{theorem}\label{thgenusunbound} Let $H$ be the elementary abelian $\ell$-group
$(\integ/\ell)^r$.  As in Situation \ref{mysit}, let $\st{\calc_n \ra \calm_n/T_n}_{n\in\nat}$ be a collection of relative
smooth proper curves of genus $g_n$ with full $\ell$-monodromy, and
let $\st{k_n}$ be a collection of finite fields, each equipped with $t_n: \spec k_n \ra T_n$.
Suppose that $\lim_{n\ra\infty}g_n=\infty$;
$\lim_{n\ra\infty}\epsilon_{\calc_n \ra \calm_n}(\ell,k_n) = 0$; and
for $n\gg 0$, $\xi(k_n) = 1$.  Then
\[
\lim_{n\ra\infty}\beta(\calc_n\ra\calm_n, k_n, \ell, H) =
\alpha(\infty,r),
\]
where $\alpha(\infty,r)$ is computed in Lemma \ref{alphainf}.
\end{theorem}

\pf The analysis is the same as that in Theorem \ref{maintheorem}.  For $n$
sufficiently large that $\xi(k_n) = 1$, we have
\[
\abs{\beta(\calc_n \ra \calm_n, k_n,\ell,H) - 
\alpha(g_n, r)} < 
\epsilon_{\calc_n \ra
    \calm_n}(\ell,k_n)
.\]
By Lemma \ref{alphainf}, $\lim_{n\ra\infty} \alpha(g_n,r)$ exists,
with limit $\alpha(\infty,r)$.  By hypothesis,
$\lim_{n\ra\infty}\epsilon_{\calc_n \ra \calm_n}(\ell,k_n)=0$; the
theorem then follows.\qed

As predicted in \cite{fw}, the divisor class groups of curves satisfy
a Cohen-Lenstra type result.  Recent research also addresses the
distribution of other ideal class groups of function fields
\cite{cardonmurty, friesen01, pacelli}.  These studies work with an
explicit affine model for a family of curves.  To ease notation
somewhat, we work with a relative curve $\calc \ra \calm/k_0$ over a
fixed finite field $k_0$, and specify an affine model by introducing a
nonempty collection of sections $S = \st{ \sigma_1, \cdots, \sigma_n :
\calm \ra \calc}$ with disjoint image.  Since we need to pass to
extension fields to apply our main result, we assume each
$\sigma_i$ is defined over the base field, $k_0$.  For a curve $C$ and a
nonempty finite set of points $S$, let $\calo_{C,S} = \cap_{P\not\in
S} \calo_P$ be the ring of functions regular outside $S$.  Let
$\cl(\calo_{C,S})$ be the ideal class group of this Dedekind domain,
and let $\cl(\calo_{C,S})_\ell$ be the $\ell$-Sylow part of that
group.

The techniques of this paper don't yield exact formulae for the
frequency with which a given group $H$ occurs as
$\cl(\calo_{C,S})_\ell$.  Still, we can at least give bounds
for the occurrence of $\ell$-Sylow subgroups of
given rank; these bounds are nontrivial if the genus of $C$ is larger than $\abs S$. Let $\rank_\ell H =
\dim_{\integ/\ell}H/\ell H$.

\begin{cor}\label{coraffine} Let $\calc\ra
\calm/k_0$ be a smooth proper relative curve of genus $g$ with full $\ell$-monodromy.
Let $S$ be a nonempty finite set of
sections $\sigma : \calm \ra \calc$ with disjoint image inside
$\calc$, and let $S_x = \cup_{\sigma\in S}\sigma(x)$.  Let $k$ be a sufficiently large finite extension of $k_0$
with $\abs{k} \equiv 1 \bmod \ell$.  For any nonnegative integer $r$,
\begin{align}
\frac{\abs{\st{x\in \calm(k) : \rank_\ell \cl(\calo_{\calc_x,S_x}) \le r}}}{\abs{\calm(k)}}
&\ge
\sum_{j = 0}^r \phi(g,j) - \epsilon_{\calc \ra \calm}(k,\ell),\label{eqsmallrank}\\
\intertext{while}
\frac{\abs{\st{x\in \calm(k) : \rank_\ell
      \cl(\calo_{\calc_x,S_x}) \ge r}}}{\abs{\calm(k)}}
& \ge \sum_{j=r+\abs S}^g \phi(g,j) - \epsilon_{\calc \ra \calm}(k,\ell).\label{eqbigrank}
\end{align}
\end{cor}

\pf Given Theorem \ref{maintheorem}, all that's necessary is to relate the
ideal class group to the divisor class group.
By a theorem of F. K. Schmidt  
\cite[proposition 1]{rosensclass}, there is an exact sequence of groups
\begin{equation}\label{classdiag}
\xymatrix{
0 \ar[r] &
0 \ar[r] & \frac{\cald(\calc_x,S_x)^0}{\calp(\calc_x,S_x)}\ar[r] & \jac(\calc_x)(k) \ar[r]&\cl(\calo_{\calc_x,S_x})\ar[r] 
& 0,
}
\end{equation}
where $\frac{\cald(\calc_x,S_x)^0}{\calp(\calc_x,S_x)}$ is the class
group of divisors of degree zero represented by divisor classes
supported at $S$.  (The sequence \eqref{classdiag} is exact on the right because the
sections $\sigma$ are defined over $k$.)  On one hand, this shows that the
$\ell$-rank of the ideal class group is no bigger than that of the the
divisor class group.  On the other hand, the $\ell$-rank of the kernel
of the surjection $\jac(\calc_x)(k) \ra \cl(\calo_{\calc_x, S_x})$ is
at most $\abs S$.  These two observations yield inequalities
\eqref{eqsmallrank} and \eqref{eqbigrank}, respectively.\qed

\section{Examples}

We conclude by working out some examples of these considerations.
Specifically, we show how Theorem \ref{maintheorem} and its variants,
in conjunction with the calculations in Section \ref{secsp}, let us
recover results of \cite{lenstra}; justify a heuristic used in
\cite{gekeler}; improve the main results of \cite{cardonmurty}, and a
special case of \cite{pacelli}; and discuss the conjecture
of \cite{fw}.

For the most part, we work over a fixed finite field $k \iso \ff_q$.
We often phrase our results in terms of $\alpha(g,r)$, which is
computed by Lemma \ref{calcphi} and Proposition \ref{calcu}.   Values of $\alpha(g,r)$
for $g \le 3$ are shown in Table \ref{bigtable}.

\begin{table}
\[
\begin{array}{||ll|l||}
\hline\hline
g & r & \alpha(g,r)\\
\hline
1 & 0 &{\frac {{\ell}^{2}-\ell-1}{{\ell}^{2}-1}}\\
1 & 1 &\frac{1}{\ell}\\
1 & 2 &{\frac {1}{\ell \left( {\ell}^{2}-1 \right) }}\\
\hline
2 & 0 &{\frac {{\ell}^{6}-{\ell}^{5}-{\ell}^{4}+\ell+1}{
 \left( {\ell}^{2}-1 \right)  \left( {\ell}^{4}-1 \right) }}\\
2 & 1 &{\frac {{\ell}^{3}-\ell-1}{{\ell}^{2} \left( {\ell}^{2}
-1 \right) }}\\
2 & 2 &{\frac {{\ell}^{3}-\ell-1}{{\ell}^{2} \left( {\ell}^{2}
-1 \right) ^{2}}}\\
2 & 3 &{\frac {1}{ \left( {\ell}^{2}-1 \right) {\ell}^{4
}}}\\
2 & 4 &{\frac {1}{{\ell}^{4} \left( {\ell}^{2}-1
 \right)  \left( {\ell}^{4}-1 \right) }}\\
\hline
3 & 0 &{\frac {{\ell}^{12}-{\ell}^{11}-{\ell}^{10}+{\ell}^{7}+
{\ell}^{5}+{\ell}^{4}-{\ell}^{3}-\ell-1}{ \left( {\ell}^{2}-1 \right)  \left( {\ell}^{4}
-1 \right)  \left( {\ell}^{6}-1 \right) }}\\
3 & 1 &{\frac {{\ell}^{8}-{\ell}^{6}+{\ell}^{2}-{\ell}^{5}+\ell-{
\ell}^{4}+1}{{\ell}^{3} \left( {\ell}^{2}-1 \right)  \left( {\ell}^{4}-1 \right) }
}\\
3 & 2 &{\frac {{\ell}^{8}-{\ell}^{6}+{\ell}^{2}-{\ell}^{5}+\ell-{
\ell}^{4}+1}{{\ell}^{3} \left( {\ell}^{2}-1 \right) ^{2} \left( {\ell}^{4}-1
 \right) }}\\
3 & 3 &{\frac {{\ell}^{5}-{\ell}^{3}-1}{{\ell}^{7} \left( {
\ell}^{2}-1 \right) ^{2}}}\\
3 & 4 &{\frac {{\ell}^{5}-{\ell}^{3}-1}{{\ell}^{7} \left( {
\ell}^{2}-1 \right) ^{2} \left( {\ell}^{4}-1 \right) }}\\
3 & 5 &{\frac {1}{ \left( {\ell}^{2}-1 \right) 
 \left( {\ell}^{4}-1 \right) {\ell}^{9}}}\\
3 & 6 &{\frac {1}{{\ell}^{9} \left( {\ell}^{2}-1
 \right)  \left( {\ell}^{4}-1 \right)  \left( {\ell}^{6}-1 \right) }}\\
\hline\hline
\end{array}
\]
\caption{\label{bigtable}The proportion of symplectic matrices of dimension $2g$ with
  fixed subspace of exact dimension $r$, as computed in Lemma
  \ref{calcphi} and Proposition \ref{calcu}.}
\end{table}

\subsection{Elliptic curves, $q\equiv 1 \bmod \ell$}

The $\ell$-torsion of a random elliptic curve $E$ is
\begin{equation}\label{eqecq1}
E[\ell](k) \iso \left\{
\begin{array}{ll}
\st 1 & \text{with probability close to } {\frac {{\ell}^{2}-\ell-1}{{\ell}^{2}-1}}\\
\integ/\ell & \text{with probability close to }\oneover\ell \\
(\integ/\ell)^2 & \text{with probability close to } \oneover{\ell(\ell^2-1)}
\end{array}
\right.
\end{equation}
in the following sense.

Let $\cale \ra \calm$ be a non-isotrivial family of elliptic curves,
such as the Legendre family (with affine model) $y^2 =
x(x-1)(x-\lambda)$ over the $\lambda$-line.  Such a family is
versal, and therefore \cite[2.2]{achhol} has full $\ell$-monodromy for
almost all $\ell$; fix one such $\ell$.  With a slight simplification
of the notation of Equation \eqref{defbeta}, let
\[
\beta_{\cale\ra\calm}(k,r) = \frac{\abs{\st{x\in\calm(k):
      \cale_x[\ell](k)\iso (\integ/\ell)^r}}}{\abs{\calm(k)}}
\]
be the proportion of elliptic curves in our family, defined over $k$,
for which the $\ell$-torsion subgroup is isomorphic to
$(\integ/\ell)^r$.  Suppose that $\abs k \equiv 1 \bmod \ell$ and
$\abs k$ is sufficiently large.  Then Theorem \ref{maintheorem} says that
\[
\abs{\beta_{\cale\ra\calm}(k,r) - \alpha(1,r)} \le
\epsilon_{\cale\ra\calm}(\ell,k),
\]
where the error term decays as $1/\sqrt{\abs k}$, and $\alpha(1,r)$,
defined in Equation \eqref{defalpha}, may be read off from the first
section of Table \ref{bigtable}.

\subsection{Elliptic curves, $q\not\equiv 1 \bmod \ell$}\label{subsececqn1}

In the situation $\cale \ra \calm$ considered above, suppose that $k$
is a large finite field for which $\abs{k}\equiv \xi \not\equiv 1
\bmod\ell$.  Again, Theorem \ref{maintheorem} says that
\begin{equation}\label{eqecqn1}
\abs{\beta_{\cale\ra\calm}(k,r) - \alpha^\xi(1,r)} \le
\epsilon_{\cale\ra\calm}(\ell,k),
\end{equation}
where
\[
\alpha^\xi(1,r) \defeq\frac{\abs{\st{x \in \gsp_2^\xi(\integ/\ell) :\ker(x-\id)
\iso (\integ/\ell)^r }}}{\abs{\sl_2(\integ/\ell)}}.
\]
We have not computed $\alpha^\xi(g,r)$ in general, but it is not hard 
to compute $\alpha^\xi(1,r)$ directly (see also \cite[3.3]{achhol}):
if $\xi\not = 1$, then 
\begin{equation}\label{eqcalcalphaxi}
\alpha^\xi(1,r) 
=
\left\{
\begin{array}{ll}
\frac{\ell-2}{\ell-1}& r = 0\\
\oneover{\ell-1}& r = 1 \\
0& r =2
\end{array}
\right..
\end{equation}
Note that this is compatible with the familiar result (use the Weil
pairing) that if $k$ has no $\ell^{th}$ root of unity, then an
elliptic curve over $k$ cannot have all its $\ell$-torsion defined
over $k$.

Taken together, Equations \eqref{eqecq1}, \eqref{eqecqn1} and
\eqref{eqcalcalphaxi}, in the special case where $k$ is a prime field
$\ff_p$, fully recover Theorem 1.14 of \cite{lenstra}.

In \cite{gekeler}, Gekeler studies the distribution of Frobenius
elements of elliptic curves over $\ff_p$, taken as elements of
$\gl_2(\integ_\ell)$.   Among
other results, he computes the proportion of elements in
$\gl_2(\integ/\ell^e)$ with given trace and determinant.  (This is
easier than the analogous question in $\sp_{2g}(\integ/\ell^e)$, first 
because conjugacy and stable conjugacy coincide in $\gl_2$, and second
because of the severe constraints on Jordan blocks of $2\times 2$
matrices.)   Combining \cite[4.4]{gekeler} and Theorem
\ref{maintheorem} allows one to compute the proportion of elliptic
curves with $\ell^e$-torsion isomorphic to a given abelian $\ell$-group
$H$.

Moreover, we can justify  a heuristic
 used in section 3 of \cite{gekeler}. There, it is asserted that if
 $m$ and $n$ are relatively prime, then for a fixed elliptic curve
 $E/k$, the actions of Frobenius on $E[m](\bar k)$ and $E[n](\bar k)$
 are independent, at least if $m\cdot n$ is small relative to
 $\sqrt{\abs k}$.  Indeed, let $\calc \ra \calm\ra k_0$ be any
 relative curve with full $mn$-monodromy; for simplicity, assume that
 $\abs{k_0}\equiv 1 \bmod mn$.  Then Frobenius elements of Jacobians
 of curves $\calc_x$ are, by Theorem \ref{katz}, equidistributed in
 $\sp_{2g}(\integ/mn)$.  Since this latter group is isomorphic to
 $\sp_{2g}(\integ/m)\oplus \sp_{2g}(\integ/n)$, Gekeler's claim
 follows.

\subsection{Quadratic function fields}\label{subseche}

%If $\deg M = d$ is even, this is called a real quadratic field.
%If $\deg M = d$ is odd, then imaginary.  In the former case, the genus is $d/2-1$; in the latter, $(d-1)/2$.  In general, the genus is $\floor{\half{d-1}}$.

Attention has recently turned to the explicit construction of ideal
 classes of given order in the class groups of quadratic function
 fields $\ff_q(x,\sqrt{f(x)})$.  Friesen computes both empirical
 \cite{friesen00} and analytic \cite{friesen01} bounds for the chance
 that $\ell$ divides the class number of $\ff_q(x)[y]/(y^2 - f(x))$,
 where $f$ is a quartic polynomial.  Cardon and Murty
 \cite{cardonmurty} show that there are at least $q^{d (\half 1 +
 \oneover \ell)}$ imaginary quadratic extensions $K = \ff_q(x,\sqrt
 {f(x)})$ of $\ff_q(x)$ where $\deg f \le d$ and the ideal class group
 of $K$ has an element of prime order $\ell \ge 3$.     While as $q$
 gets large this produces arbitrarily large families of quadratic
 function fields with class number divisible by $\ell$, it is a
 vanishingly small proportion of all quadratic function fields.

We can use Corollary \ref{coraffine} to compute the {\em
 proportion} of quadratic function fields with class number divisible
 by $\ell$, and thereby strengthen these results.

Suppose $q$ is a power of an odd prime and that $q\equiv 1 \bmod
 \ell$.  We let $k = \ff_q$, and let $\st{k_n}$ be any collection of
 finite extensions of $k$ with $\lim_{n\ra\infty} \abs{k_n} = \infty$.

Let $\calh_d$ be the space of separable monic polynomials $f(x)$ of
degree $d$.  Over it lies $\calc_d$, the curve with affine model $y^2
= f(x)$; it is a hyperelliptic curve of genus $\floor{\half{d-1}}$.
We work under the hypothesis that $\calc_d \ra \calh_d$ has full
$\ell$-monodromy.  For
odd $d$, this is implied by unpublished work of J.K. Yu
\cite[10.5.10]{katzsarnak};  we will treat the general case in a
future work. 

The function field of the curve with affine model $y^2 = f(x)$ is
called an imaginary quadratic function field if $d = \deg f$ is odd,
and a real quadratic function field otherwise.  We address these cases
separately.

If $d$ is odd, then there is a single point ``at infinity'' in this
affine model; the left-hand term of \ref{classdiag} is
trivial, and the ideal class group of this ring is isomorphic to the
$\ff_q$-rational points of the Jacobian of the associated proper
curve.  We see that, for instance,
\[
 \lim_{n\ra\infty}\frac{\abs{\st{f(x) \in k_n[x]: \deg f = d,
       f\text{ monic}, \ell |
       \cl(k_n[x,\sqrt{f(x)}])}}}
 {{\abs{\st{f(x) \in k_n[x] :
       \deg f = d, f\text{ monic}}}}}
\]
is equal to
\[
\lim_{n\ra\infty}\frac{\abs{\st{f(x) \in k_n[x]: \deg f = d,
       f\text{ monic, separable}, \ell |
        \cl(k_n[x,\sqrt{f(x)}])}}}
  {{\abs{\st{f(x) \in k_n[x] :
        \deg f = d, f\text{ monic, separable}}}}},
\]
since most polynomials are separable, which is in turn equal to 
\[
\lim_{n\ra\infty} \frac{\abs{ \st{ x\in \calh_d(k_n): \ell |
        \abs{\jac(\calc_x)[\ell](k_n)}}}}{\abs{\calh_d(k_n)}},
\]
or $1-\alpha(g,0)$.

If $d$ is even, then there are two points at infinity, and the
regulator term in \ref{classdiag} is an abelian group on a single
generator.  
Therefore, for any curve $C/k_n$ with affine model $C^{{\rm aff}}: y^2 =
f(x)$, we have
\[
\rank_\ell(\jac(C)[\ell](k_n)) \ge
\rank_\ell(\cl(\calo_{C^{{\rm aff}}}))\ge
\rank_\ell(\jac(C)[\ell](k_n)) - 1,
\]
and the chance that $\ell$ divides the class group of the affine
coordinate ring is bounded from below by 
\[
\sum_{r=2}^g \alpha(g,r) - \epsilon_{\calc_d \ra \calh_d}(k_n,\ell).
\]
Thus, we can significantly strengthen the main conclusion of
\cite{cardonmurty}; as $n\ra\infty$, there is an element of order
$\ell$ in the class group of $k_n[x][y]/(y^2 -f(x))$ for a positive
proportion of monic degree $d$ polynomials $f(x)\in k_n[x]$.

Note that Theorem \ref{thgenusunbound} allows one to make uniform statements about
curves of the form $y^2 = f(x)$ as $\deg f \ra \infty$, provided that
the size of the base field grows sufficiently quickly.   

\subsection{Cyclic cubic fields}\label{extrigonal}
Pacelli \cite{pacelli} looks at curves $y^d = f(x)$, and obtains
results (for general $d$) similar to those of \cite{cardonmurty} for $d=2$.
As mentioned in the introduction, the general family of curves (with
affine model) $y^d = f(x)$ cannot have full $\bmod\ \ell$ monodromy,
because of the extra automorphisms this family possesses.  Still, by
computing in the appropriate monodromy group one can calculate
divisibility of class numbers for these families.  We expect
\cite{zarendjac} that the monodromy group is a unitary group
associated to $\rat(\zeta_d)$.

We take up these considerations in the special case where $d=3$, the
degree of $f$ is $4$, and $3$ is invertible in the base field.  Let
$\calc \ra \calp$ be the family of curves with affine model $y^3 =
f(x)$, where $f$ ranges over all separable polynomials of degree $4$.
Each fiber $\calc_x$ has genus $3$.  Moreover, since there is an
obvious action of a cyclic group of order $3$ on $\calc$, the Jacobian
$\jac(\calc)$ admits an action by $\integ[\zeta_3]$.  The action on
the tangent space at the identity of any fiber has signature $(2,1)$,
since actions of type $(3,0)$ are rigid.  Therefore, under the Torelli
map, $\calc \ra\calp$ becomes identified with an open subset of the
Picard modular variety associated to $\integ[\zeta_3]$.  Using
transcendental arguments \cite{holzapfel} and the theory of
compactification \cite {larsenthesis}, one knows that for almost all
$\ell$, the full $\ell$-adic monodromy group of this family is
$G(\integ_\ell)$, where $G$ is the unitary group in three variables
associated to $\integ[\zeta_3]$.

Suppose, then, that $\calc \ra \calp$ has $\ell$-monodromy group
$G(\integ/\ell)$.  If $\integ[\zeta_3]$ is inert at $\ell$, then
$G(\integ/\ell)$ is an example of the unitary groups studied in
Section \ref{subsecunitary}.  In particular, we see that:
\[
\jac(y^3 =f(x))[\ell](k) \iso
 \left\{
\begin{array}{ll}
\st 1 & \text{with probability close to } {\frac {\ell \left( {\ell}^{5}-{\ell}^{3}-1 \right) 
}{ \left( \ell+1 \right)  \left( {\ell}^{2}-1 \right)  \left( {\ell}^{3}+1
 \right) }}\\
(\integ/\ell)^2 & \text{with probability close to } {\frac {{\ell}^{5}-{\ell}^{2}+{\ell}^{4}-\ell-1}{{\ell}^{2
} \left( \ell+1 \right) ^{2} \left( {\ell}^{2}-1 \right) }}\\
(\integ/\ell)^4 & \text{with probability close to }{\frac {{\ell}^{3}+{\ell}^{2}-1}{ \left( \ell+1
 \right) ^{2} \left( {\ell}^{2}-1 \right) {\ell}^{3}}}\\
(\integ/\ell)^6 & \text{with probability close to }{\frac {1}{{\ell}^{3} \left( \ell+1 \right) 
 \left( {\ell}^{2}-1 \right)  \left( {\ell}^{3}+1 \right) }}
\end{array}
\right..
\]
(If $\integ[\zeta_d]$ splits at $\ell$, then $G(\integ_\ell)$ is
isomorphic to a general linear group, and a similar, but easier,
calculation applies.)

\subsection{The Friedman-Washington conjecture}\label{subsecfw}

Let $\calc \ra \calm$ be a family of curves of genus $g$ with full
$\ell$-monodromy; this corresponds to any suitably general family of
curves.  Let $k$ be a large finite field, say with $\abs k \equiv 1
\bmod \ell$.  Then Theorem \ref{maintheorem} says that the proportion
of fibers $\calc_x$, for $x\in \calm(k)$, with $\jac(\calc_x)[\ell](k)
\iso (\integ/\ell)^r$ is $\alpha(g,r)$; see Table \ref{bigtable} for
the first few values of $\alpha(g,r)$.

Friedman and Washington \cite{fw} give a conjectural
description of the frequency with which a given abelian $\ell$-group
occurs as the $\ell$-Sylow part of the divisor class group of a
function field.  While they formulate their conjecture in terms of
hyperelliptic curves in order to preserve the analogy with the Cohen-Lenstra
heuristics, all of their arguments depend on merely having a
sufficiently general family of curves.  Given our
expectation (Section \ref{subseche}) that hyperelliptic curves behave,
in terms of $\ell$-monodromy, like general curves,  we compare the
predictions of \cite{fw} to the results of Theorem \ref{maintheorem}
for $\calc\ra\calm$ with full monodromy.

To facilitate this comparison we estimate the chance that the
$\ell$-part of the class group of a curve is trivial.  Let
$\phi_{\gl}(n)$ denote the proportion of elements $x\in \gl_n$ for
which $x-\id$ is invertible.  It is shown that for large $n$ $\phi_{\gl}(n)$
approaches
\[
\til\phi_{\gl}(n) \defeq \prod_{j=1}^n ( 1 - \ell^{-j}).
\]
In the special case of genus $2$, Friedman and Washington
predict that the proportion of curves with trivial $\ell$-class group
is (close to) $\til \phi_{\gl}(4)$, while Theorem \ref{maintheorem} says
that this proportion is actually $\alpha(2,0) = {\frac {{\ell}^{6}-{\ell}^{5}-{\ell}^{4}+\ell+1}{
 \left( {\ell}^{2}-1 \right)  \left( {\ell}^{4}-1 \right) }}$.  The
gap between $\til\phi_{\gl}(4)$ and $\alpha(2,0)$ is of order
$1/\ell^2$.

Now, \cite[p.131]{fw} expresses the hope that this discrepancy disappears for large
genus;  unfortunately, this difference persists.  The proportion of curves of
(arbitrarily large) genus with trivial $\ell$-group approaches (the
well-defined limit; see Lemma \ref{alphainf}) $\phi(\infty)$, which by
\cite[3.3]{achhol} is $1 -\ell/(\ell^2-1) + O(1/\ell^3)$.  The difference
between the conjectural estimate of \cite{fw} and the actual value
remains of order $1/\ell^2$, even as the genus of the curves in
question gets arbitrarily large.

\bibliographystyle{plain}
\bibliography{jda}

\end{document}